\documentclass[12pt]{amsart}
\usepackage{amssymb}
\usepackage{amscd}
\usepackage{a4}

\newcounter{liste}

\newenvironment{lister}{\begin{list}%
{{\rm(\roman{liste}) }}{%
\usecounter{liste}%
\setlength{\labelsep}{0ex}%
\settowidth{\labelwidth}{(viii)\ }%
\settowidth{\leftmargin}{(viii)\ }%
\addtolength{\topsep}{1\topsep}%
}}{\end{list}}

\newenvironment{listerp}{\begin{list}%
{{\rm(\roman{liste})' }}{%
\usecounter{liste}%
\setlength{\labelsep}{0ex}%
\settowidth{\labelwidth}{(viii)'\ }%
\settowidth{\leftmargin}{(viii)'\ }%
\addtolength{\topsep}{1\topsep}%
}}{\end{list}}

\begin{document}

\noindent This work appeared in ``Continuum theory: In Honor of Professor David P.\ Bellamy on the Occasion of His 60th Birthday'',
I.\ W.\ Lewis, S.~Mac\'\i as, S.\ B.\ Nadler (eds.), Aportaciones Mat.\ Investig.\ {\bf 19} (2007), pp.\ 71--84.

\vspace{1.1cm}

\title[The non-existence of  common models]{The non-existence of  common models for some classes of
 higher-dimensional hereditarily indecomposable continua}

\author[Jerzy Krzempek and El\.{z}bieta Pol]{Jerzy Krzempek$^1$ and El\.{z}bieta Pol$^1$}

\thanks{$^1$Research partially supported by MNiSW Grant Nr. N201 034 31/2717.}

\thanks{Dedicated to D.P.Bellamy on the occasion of his 60th birthday}

\address{Institute of Mathematics\\
Silesian University of Technology\\
Kaszubska 23\\
44-100 Gliwice\\
Poland}

\email{jkrzempek.math@gmail.com}

\address{Institute of Mathematics\\
University of Warsaw\\
Banacha 2\\
02-197 Warszawa\\
Poland}

\email{pol@mimuw.edu.pl}

\begin{abstract}
A continuum $K$ is a common model for the family ${\mathcal K}$ of continua if every member of ${\mathcal K}$ is a continuous image of
$K$. We show that none of the following classes of spaces has a common model: 1) the class of strongly chaotic hereditarily
indecompos\-able $n$-dimen\-sion\-al Cantor manifolds, for any given natural number $n$, 2) the class of strong\-ly chaotic
hereditarily indecomposable hereditarily strongly infinite-dimensional Cantor manifolds, 3) the class of strongly chaotic hereditarily
indecomposable continua with\linebreak transfinite dimension (small or large) equal to $\alpha$, for any given ordinal number $\alpha
< \omega_{1}$.
\end{abstract}




\maketitle

\section{Introduction}
Our terminology follows \cite{En}. By {\it dimension} we
understand the covering dimension dim unless otherwise stated. By
a {\it continuum} we mean a compact connected space (we assume
that all our spaces are metrizable and separable). A continuum $X$
is {\it indecomposable}, if it is not the union of two
 proper subcontinua. A continuum $X$ is {\it hereditarily
indecomposable}, abbreviated  HI, if every subcontinuum of $X$ is
indecomposable. A continuum $X$ is {\it strongly chaotic}, if for
any two disjoint subsets $U$ and $V$ of $X$, with $U$ being
nonempty and open, there is no homeomorphism from $U$ onto $V$.
All our mappings are meant to be continuous.

We will say that a continuum $K$ is a {\it common model} for the
class ${\mathcal K }$ of continua, if each member of ${\mathcal K
}$ is a continuous image of $K$ (we do not assume that $K \in
{\mathcal K }$). Z.Waraszkiewicz \cite{Wa} constructed in 1932 a
family of plane continua (called the Waraszkiewicz spirals)
without a common model. Applying this result, D.P.Bellamy
\cite{Be} has shown that the collection of all indecomposable
continua has no common model. Other collections of indecomposable
continua without a common model were constructed, among others, by
R.L.Russo \cite{Rus} (collection of planar indecomposable
tree-like continua), W.T.Ingram \cite{In} (collection of planar
hereditarily indecomposable tree-like continua) and
T.Ma\'{c}kowiak and E.D.Tymchatyn \cite{MT} (collection of
1-dimensional hereditarily indecomposable continua). More
information concerning classes of continua with or without a
common model can be found in \cite{MT}.

 In this paper we will show that there is no
common model for the class ${\mathcal K}$ of strongly chaotic
hereditarily indecomposable continua which are either a)
$n$-dimensional Cantor manifolds, for any given natural number
$n$, or  b) here\-di\-ta\-rily strongly infinite-dimensional
Cantor manifolds, or else c) countable-dimensional continua with
transfinite dimension (small or large) equal to $\alpha$, for any
given $\alpha < \omega _{1} $.

 For any of the collections ${\mathcal
K}$ we have described, we will find a family of Waraszkiewicz
spirals without a common model such that for every spiral $W$ from
this family there  exists a continuum in ${\mathcal K}$ which can
be mapped (continuously) onto $W$.

To this end, we shall apply, among others, some results from
\cite{PE3} and a theorem of J.T.Rogers, Jr. from \cite{Ro}.

In the last section we will strengthen some of the results from
the preceding section using
  a method suggested by an anonymous referee of the
paper \cite{PE1} (see Remark 5.2 in \cite{PE1}). Namely, we will show that for every collection ${\mathcal K}$ mentioned
above and for every compactification $W$ of the ray there exists a continuum in ${\mathcal K}$ which can be mapped onto~$W$.

\section{Preliminaries} The first hereditarily indecomposable continuum, now called the pseudoarc, was
constructed by B.Knaster \cite{Kn}. The first examples of
$n$-dimensional HI continua, for every $n = 2,3,\ldots ,\infty $,
were given by R.H.Bing \cite{Bi1} in 1951. Several constructions
of collections of non-homeomorphic higher-dimensional HI continua,
with different additional properties, were given in \cite{PE1},
 \cite{PE3}, \cite{PE4}, \cite{PR} and \cite{RM}.

 A compact $n$-dimensional space is an {\it $n$-dimensional Cantor
 manifold} if no closed subset $L$ of $X$ with dim$L \leq n-2$
 disconnects $X$ (see \cite{En}, sec.1.9).

Below by the {\it ray} we will understand a space  homeomorphic to
the half-line $[0,+\infty )$. By a {\it Waraszkiewicz spiral} we
mean  a member of a family of planar continua without a common
model, constructed by Z.Waraszkiewicz \cite{Wa}. Every
Waraszkiewicz spiral is a compactification of the ray with the
remainder homeomorphic to the circle.

\medskip

{\bf Definition 2.1.} Below by $\{W_{K}: K \in {\mathcal S}\}$, where ${\mathcal S}$ is some set of cardinality $2^{\aleph
_{0}}$, we will denote a family consisting of some Waraszkiewicz spirals such that \medskip

(1) no continuum $X$ can be mapped onto uncountably many of the
spirals $W_{K}$.

\medskip The proof that there exists such a family, can be found, for example, in the paper of T.Ma\'{c}kowiak and
E.D.Tym\-cha\-tyn \cite{MT} (see (20.3), (20.4) and (20.9)). As pointed out by the authors, see \cite{MT}, page 49, their
proof is based on an unpublished construction by D.P.Bellamy. A new construction of a family of Waraszkiewicz spirals
without a common model, satisfying the condition (1), was also presented by D.P.Bellamy at the 2004 Spring Topology nad
Dynamics Conference in Birmingham, Alabama.

Let us recall also the following theorem of J.T.Rogers, Jr. (see
\cite{Ro}, Theorem 4):

\medskip

 (2) for each Waraszkiewicz spiral $W$ there exists  a 1-dimensional plane continuum
$\hat{W}$ such that each nondegenerate subcontinuum  of $\hat{W}$
can be mapped onto $W$.

\medskip

Using this theorem, M.Re\'{n}ska \cite{RM} has proved that

\medskip

(3)  for every natural number $n$ and for every Waraszkiewicz
spiral $W$ there exists an HI $n$-dimensional Cantor manifold
every $n$-dimensional subcontinuum of which can be mapped onto
$W$.

\medskip

Although this fact was not stated explicitly by Re\'{n}ska, its
proof can be easily extracted from the proof of Lemma 3.4 in
\cite{RM}. From (3) it follows that there is no common model for
HI $n$-dimensional Cantor manifolds.

Using  the result (2) of Rogers in a different way, we will obtain
the following strengthening of (3).

\medskip

{\bf  Lemma 2.2.}   {\it For every natural number $n$ and for
every Waraszkiewicz spiral $W$ there exists an HI $n$-dimensional
Cantor manifold $X(W)$ such that every nontrivial subcontinuum of
$X(W)$ can be mapped onto $W$.}

\medskip

{\bf Proof.} Fix a Waraszkiewicz spiral $W$ and let $\hat{W}$ be a
continuum satisfying (2), constructed by Rogers. Let
  $X = \prod _{i=1}^{n+1}X_{i}$, where
$X_{i}=\hat{W}$ for every $i$, and let $p_{i}: X \to X_{i}$ be a
projection. Then $X$ is $(n+1)$-dimensional (see \cite{En},
Problem 1.8.K (b)), so, by a theorem of Bing \cite{Bi1}, it
contains an $n$-dimensional HI continuum $Y$. Let $X(W)\subset Y $
be an $n$-dimensional Cantor manifold (see \cite{En}, Theorem
1.9.9). Then for any nontrivial subcontinuum $Z$ of $X(W)$, the
projection $p_{i}(Z)$  is nontrivial for some $i$. Since the
continuum $p_{i}(Z)$ can be mapped onto $W$, the same is true for
$Z$.

\medskip

We shall recall now some notions concerning infinite-dimensional
spaces.  A space $X$ is {\it strongly infinite-dimensional}, if
there exists an infinite sequence $(A_{1},B_{1}), (A_{2},B_{2}),
\ldots $ of pairs of disjoint closed sets in $X$, such that if
$L_{i}$ is a partition between $A_{i}$ and $B_{i}$ in $X$, then
$\bigcap _{i=1}^{\infty}L_{i}\neq \emptyset$. A space $X$ is {\it
hereditarily strongly infinite-dimensional},  if every subspace of
$X$  is either 0-dimensional or strongly infinite-dimensional. The
first example of a hereditarily strongly infinite-dimensional
continuum was given by L.Rubin \cite{Rub} (cf. also  \cite{En},
Problem 6.1.C). An infinite-dimensional continuum $X$ is an {\it
infinite-dimensional Cantor manifold} if every closed subset of
$X$ which disconnects $X$ is infinite-dimensional. An
infinite-dimensional Cantor manifold, which is a hereditarily
strongly infinite-dimensional continuum, will be called a
hereditarily strongly infinite-dimensional Cantor manifold. The
{\it small (large) transfinite dimension} ind (Ind) is the
transfinite extension of the classical small (large) inductive
dimension (see \cite{En}, Definitions 7.1.1 and 7.1.11). A compact
metrizable space $X$ has small (equivalently, large) transfinite
dimension if and only if it is countable dimensional, i.e.\ it is
the union of countably many finite-dimensional sets.

\medskip

{\bf  Lemma 2.3.}  {\it For every Waraszkiewicz spiral $W$ there
exists an HI here\-di\-tarily strongly infinite-dimensional Cantor
manifold $X(W)$ such that every nontrivial subcontinuum of $X(W)$
can be mapped onto $W$.}

\medskip

{\bf Proof.} Fix a spiral $W$. Let $\hat{W}$ be a continuum
satisfying (2), constructed by Rogers. Then the countable product
$X = \prod _{i=1}^{\infty}X_{i}$, where $X_{i}=\hat{W}$ for every
$i$, is strongly infinite-dimensional (see \cite{Li} or \cite{En},
Problem 6.1.H (b)). By a theorem of Rubin \cite{Rub}, $X$ contains
a hereditarily strongly infinite-dimensional compactum $Y$ (see
also \cite{En}, Problem 6.1.G). By a theorem of Tumarkin
\cite{Tu}, every hereditarily strongly infinite-dimensional
compactum contains an  infinite-dimensional Cantor manifold and by
a theorem of Bing \cite{Bi1}, every continuum  can be separated by
a closed set all of whose components are hereditarily
indecomposable. Thus $Y$ contains an HI hereditarily strongly
infinite-dimensional Cantor manifold $Z$. Since $Z \subset \prod
_{i=1}^{\infty}X_{i}$, then every nontrivial subcontinuum of $Z$
has a nontrivial projection onto some $X_{i}$, hence it can be
mapped onto $W$.

\medskip

{\bf  Remark 2.4.} Let $\{W_{K}: K \in {\mathcal S}\}$  be a family of Waraszkiewicz spirals from Definition 2.1 and let $n$ be any
natural number. For $K \in {\mathcal S}$,  let $X(W_{K})$ be an HI $n$-dimensional Cantor manifold constructed in Lemma 2.2 for $W =
W_{K}$ (respectively, let $X(W_{K})$ be an HI hereditarily strongly infinite-dimensional Cantor manifold constructed in Lemma 2.3 for
$W = W_{K}$ ). Then

\medskip

(4) no continuum $X$ has  non-constant mappings into uncountably
many of the continua $X(W_{K})$.

\medskip

Thus,  for every countable family $\{X_{i}: i = 1,2,\ldots \}$ of
continua, there exists an HI $n$-dimensional Cantor manifold
(respectively, an HI hereditarily strongly infinite-dimensional
Cantor manifold) $X(W_{K})$, where $K$ is some element of
${\mathcal S} $, such that no $X_{i}$ can be mapped onto a
nontrivial subcontinuum of $X(W_{K})$.

 Applying the transfinite induction, we obtain that

 \medskip

 (5)  there exists a family
$\{X_{\alpha}: \alpha < \omega _{1}\}$ of HI $n$-dimensional
Cantor ma\-ni\-folds (respectively, of  HI hereditarily strongly
infinite-dimensional Cantor manifolds) such that for every $\alpha
< \beta < \omega _{1}$, every mapping from $ X_{\alpha}$ to $
X_{\beta}$ is constant.

\medskip

We shall present now a method of condensation of singularities,
 which will be used in our constructions.

 The method of  condensation of singularities
 goes back to
Z.Janiszewski \cite{J}. It was refined by R.D.Anderson and
G.Choquet \cite{AC} and developed essentially by T.Ma\'{c}kowiak
\cite{Ma1}, \cite{Ma2}. It is closely related to V.Fedorchuk's
``method of resolutions'' \cite{Fe}, modified also by V.A.Chatyrko
\cite{Ch}.

A subcontinuum $Y$ of a continuum $X$ is {\it terminal}, if every
subcontinuum of $X$ which intersects both $Y$ and its complement
must contain $Y$. A continuum $X$ is HI if and only if every
subcontinuum of $X$ is terminal in $X$. A mapping $p : X \to Y $
is {\it atomic}, if every fiber $p^{-1}(y)$ is a terminal
subcontinuum of $X$.

A {\it composant} of a point $x$ in a continuum $X$ is the union
of all proper subcontinua of $X$ containing $x$. Composants are
dense connected $F_{\sigma}$-subsets of $X$ and every HI continuum
$X$ has $2^{\aleph _{0}}$ composants which are pairwise disjoint
and boundary in $X$.

\medskip

{\bf  Theorem 2.5.} {\it Let $X$ be a continuum,  $\{K_{i}:i \in N\}$ a sequence of continua and  $\{a_{i}: i \in N\}$ a
sequence of different  points of $X$. Then there exist a continuum $L(X,K_{i},a_{i})$ and an atomic mapping $p:
L(X,K_{i},a_{i}) \to X$ such that if $A = \bigcup _{i=1}^{\infty}\{ a_{i}\}$, then
\begin{lister}
\item $p^{-1}(a_{i})$ is a boundary set homeomorphic to $K_{i}$, for $i \in N$,
\item $p \mid p^{-1}(X \setminus A): p^{-1}(X \setminus A) \to  X \setminus A$ is a homeomorphism,
\item if $A$ is dense in $X$ then every non-empty open subset of $L(X,K_{i},a_{i})$ contains $p^{-1}(a_{i})$ for some $i \in N$,
\item if  $X$ and all $K_{i}$ are HI then $L(X,K_{i},a_{i})$ is HI,
\item if $n$ and $\alpha$ are ordinal numbers such that ${\rm ind}X \leq n < \omega _{0}$ and $n \leq {\rm ind}K_{i}
\leq \alpha < \omega _{1}$ for every $i \in N$ then ${\rm ind}L(X,K_{i},a_{i}) \leq \alpha$, and the same is true if we
replace {\rm ind} by {\rm Ind,}
\item if $X$ is an $n$-dimensional Cantor manifold, ${\rm dim}K_{i} \leq n$ and for every composant $L$ of $X$ the set
$L \cap A$ is finite, then $L(X,K_{i},a_{i})$ is an $n$-dimensional Cantor manifold,
\item if $X$ and all $K_{i}$ are hereditarily strongly infinite-dimensional, then so is $L(X,K_{i},a_{i})$,
\item if $X$ is an infinite-dimensional Cantor ma\-ni\-fold then so is $L(X,K_{i},a_{i})$.
\end{lister}}

\smallskip

For the proof of this theorem see \cite{PR}, Theorem 3.2 and
\cite{PE3}, Theorem 2.1. Roughly speaking, the space
$L(X,K_{i},a_{i})$ is obtained from continuum $X$ by replacing
points $a_{i}$ by continua $K_{i}$.

\medskip

{\bf  Remark 2.6.} In the case when $a_{1}=a$, $K_{1}=A$ and all
continua $K_{i}$ for $i \geq 2$ are degenerate, we will denote the
space $L(X,K_{i},a_{i})$ by $M(X,A,a)$. More generally, given two
continua $X$ and $A$ and a point $a \in X$, we will denote by
$M(X,A,a)$ and call a {\it pseudosuspension of $A$ over $X$ at a
point $a$} (cf. \cite{Ma1}) every continuum $M(X,A,a)$ such that
there exists an atomic mapping $p : M(X,A,a) \to X$ onto $X$ such
that $p^{-1}(a)$ is homeomorphic to $A$ and $p \mid p^{-1}(X
\setminus \{a\} ) : p^{-1}(X \setminus \{a\} ) \to X \setminus
\{a\}$ is a homeomorphism.

Using the fact, that $p^{-1}(a)$ is a terminal continuum in $M(X,A,a)$, one can easily show that
\begin{listerp}\setcounter{liste}{3}
\item if  $X$ and $A$ are HI then  so is $M(X,A,a)$,
\item if $n$ and $\alpha$ are ordinal numbers such that ${\rm ind}X \leq n < \omega _{0}$ and $n \leq {\rm ind}A \leq \alpha <
\omega _{1}$ then ${\rm ind}M(X,A,a) \leq \alpha$, and the same is true if we replace ind by Ind,
\item if $X$ is an $n$-dimensional Cantor manifold and   dim$A \leq n$ then $M(X,A,a)$ is an $n$-dimensional Cantor manifold,
\item if $X$ and  $A$ are hereditarily strongly infinite-dimensional, then so is
$M(X,A,a)$,
\item if $X$ is an infinite-dimensional Cantor ma\-ni\-fold then so is $M(X,A,a)$.
\end{listerp}

\section{Strongly chaotic HI continua}

Strongly chaotic spaces described in the Introduction were defined
by J.J.Charatonik and W.J.Charatonik in \cite{CC}.  A space $X$ is
{\it strongly rigid}, if the only embedding of $X$ into itself is
the identity of $X$ onto $X$. Every strongly chaotic space is
strongly rigid (see \cite{CC}). In \cite{Co} H.Cook proved that

\medskip

(6) there exists a one-dimensional HI continuum $E$ such that for
any two different nontrivial subcontinua of $E$, there is no
mapping from one onto the other.

\medskip

In particular, taking a family of $2^{\aleph_{0}}$ disjoint
subcontinua of the Cook continuum $E$ one obtains a family of
1-dimensional strongly chaotic HI continua such that every mapping
between distinct continua that belong to this family is constant.
Concerning higher-dimensional HI continua, the following result
was essentially proved  in \cite{PE3}.

\medskip

{\bf Theorem 3.1.} {\it  There exists a fa\-mi\-ly $\{Y_{\sigma}: \sigma \in S\}$, where $\mid S \mid = 2^{\aleph _{0}}$,
consisting of HI continua which are either
\begin{lister}
\item[\em (a) ] $n$-dimensional Cantor manifolds, for any given natural number $n$, or
\item[\em (b) ] hereditarily strongly infinite-dimensional Cantor manifolds, or else
\item[\em (c) ] continua with small (large) transfinite dimension equal to $\alpha$ every open subset of which is
infinite-dimensional, for any given infinite ordinal $ \alpha <\omega _{1}$,
\end{lister}
such that for every $\sigma , \tau\in S$,
\begin{lister}
\item every continuum $Y_{\sigma}$ is strongly chaotic,
\item no open nonempty subset of $Y_{\sigma}$ embeds in $Y_{\tau}$ if $\sigma \neq
 \tau$.
\end{lister}}

 \smallskip

 More precisely, the paper \cite{PE3} contains constructions of such families, but instead of the condition (i)  a
 weaker condition  is proved that
 every $Y_{\sigma}$ is strongly rigid (the fact that the constructed spaces are also
 strongly chaotic is stated in
 Added in proof). For the sake of completeness, we
 shall add to the results from \cite{PE3} some arguments which
 simplify the original reasoning, and yield the stronger
 conclusion we shall need.

The following proposition  follows from Theorems 1.1, 1.2, 1.3 of
\cite{PE3}, but a simpler proof can be obtained directly from
Lemmas 3.1, 4.3, 5.2 and 6.1 in \cite{PE3}.

\medskip

{\bf  Proposition 3.2.} {\it Let ${\mathcal L}$ be one of the following  classes of spaces:
\begin{lister}
\item[\em 1) ] the class ${\mathcal L}_{n}$ of $n$-dimensional HI Cantor ma\-ni\-folds, where $n=1,2,\ldots$,
\item[\em 2) ] the class ${\mathcal L}_{\infty}$ of HI hereditarily strongly infinite-dimensional  Cantor ma\-ni\-folds,
\item[\em 3) ] the class ${\mathcal L}_{\alpha}({\rm ind})$
(${\mathcal L}_{\alpha}({\rm Ind})$) of countable-dimensional HI continua with small (large) transfinite dimension equal to
$\alpha$, where $\omega _{0}\leq \alpha < \omega _{1} $.
\end{lister}

Then there exists a family $\{X_{t}: t \in T\}$ of cardinality
$2^{\aleph _{0}}$ consisting of continua from ${\mathcal L}$ such
that if $t \neq t'$ then $X_{t}$ does not embed in $X_{t'}$.}

\medskip

{\bf  Proposition 3.3.} {\it  Let $\{ X_{t}: t \in T \}$ be an
infinite family of  continua such that $X_{t}$ does not embed in
$X_{t'}$ for any $t \neq t' $.
 Let $X$ be any  continuum such that no $X_{t}$ embeds in
$X$. Let $a_{1},a_{2},\ldots$ be a sequence of points  such that the set $A = \bigcup _{i=1}^{\infty}\{a_{i}\}$ is dense in
$X$. Then
\begin{lister}
\item for every sequence $\sigma =\{s_{i}\}_{i=1}^{\infty}$ of different elements of $T$, the   continuum
$Y_{\sigma}=L(X,X_{s_{i}},a_{i})$ is strongly chaotic,
\item if $\sigma = \{s_{i}\}_{i=1}^{\infty}$ and $ \tau = \{t_{i}\}_{i=1}^{\infty}$ are disjoint sequences of different
elements of $T$ then no open nonempty subset of $Y_{\sigma}$ embeds into $Y_{\tau}$.
\end{lister}}

\smallskip

{\bf Proof.} Let $\sigma = \{s_{i}\}_{i=1}^{\infty}$ and $ \tau =
\{t_{i}\}_{i=1}^{\infty}$ be disjoint sequences of different
elements of $S$. Note that $Y_{\sigma}$ is the union of a subset
homeomorphic to a subspace of $X$ and topological copies
$X_{s_{i}}'$ of $X_{s_{i}}$ for $i=1,2,\ldots $. Since $A $ is
dense in $X$ then every nonempty open subset of $Y_{\sigma}$
contains some $X_{s_{i}}'$.

 Let $U$ be a nonempty open subset of $Y_{\sigma}$ and $h: U \to
Y_{\sigma} \setminus U$ be an embedding. Then $X_{s_{i}}' \subset
U$ for some $i$. Since $X_{s_{i}}'$ does not embed in $X$ then
$h(X_{s_{i}}')$ intersects a copy $X_{s_{j}}'$ of some
$X_{s_{j}}$, where $j \neq i$. Since $X_{s_{j}}'$ is a terminal
continuum in $Y_{\sigma}$, then $h(X_{s_{i}}') \subset X_{s_{j}}'$
or $h(X_{s_{i}}') \supset X_{s_{j}}'$ - a contradiction.

  Let $U$ be a nonempty open subset of $Y_{\sigma}$ and $h : U \to
Y_{\tau}$ be an embedding. Then  $X_{s_{i}}' \subset U$ for some
$i$. Since $X_{s_{i}}'$ does not embed in $X$ then $h(X_{s_{i}}')$
intersects a copy $X_{t_{j}}'$ of some $X_{t_{j}}$. Since
$X_{t_{j}}'$ is a terminal continuum in $Y_{\tau}$, then
$h(X_{s_{i}}') \subset X_{t_{j}}'$ or vice versa - a
contradiction.

\medskip

{\bf Proof of Theorem 3.1. } Let ${\mathcal L}$  and  $\{ X_{t}: t
\in T \}$ be as in Proposition 3.2.  Choose $t_{0}\in T$ and let
$S$ be a family of cardinality $2^{\aleph _{0}}$ consisting of
sequences $\sigma = \{s_{i}\}_{i=1}^{\infty}$  of different
elements of $T \setminus \{t_{0}\}$ such that for $\sigma =
\{s_{i}\}_{i=1}^{\infty} \neq \tau = \{t_{i}\}_{i=1}^{\infty}$ the
sequences $\{s_{i}\}_{i=1}^{\infty}$ and
$\{t_{i}\}_{i=1}^{\infty}$ are disjoint.  In the case of
${\mathcal L }={\mathcal L}_{n}$ or ${\mathcal L }={\mathcal
L}_{\infty}$ let $X = X_{t_{0}}$. Additionally, in the case of
${\mathcal L}_{n}$ we assume that the points $a_{1},a_{2}, \ldots$
belong to different composants of $X$. Then the family
$\{Y_{\sigma}: \sigma \in S\}$ consisting of continua $Y_{\sigma}$
constructed in Proposition 3.3 satisfies the conditions of Theorem
3.1, since by Theorem 2.5 (iv), (vi), (vii) and (viii), if all
$X_{t}$ and $X$ belong to ${\mathcal L}$ then $Y_{\sigma} \in
{\mathcal L}$. In the case when ${\mathcal L }={\mathcal
L}_{\alpha}({\rm ind})$ (or ${\mathcal L }={\mathcal
L}_{\alpha}({\rm Ind})$), let $X$ be any 1-dimensional HI
continuum. Then the family $\{Y_{\sigma}: \sigma \in S\}$
constructed in Proposition 3.3 satisfies the required conditions,
since by Theorem 2.5 (iii), (iv) and (v), if all $X_{t}$ belong to
${\mathcal L}$ and $X$ is finite-dimensional then every
$Y_{\sigma}$ belongs to ${\mathcal L}$ and every open subset of
$Y_{\sigma}$ is infinite-dimensional, because it contains a copy
of some $X_{t}$.

\medskip

{\bf  Remark 3.4.} (A) As observed in Added in Proof in \cite{PE3}
(where Lemma 3.1 should be changed to Lemma 3.2),  the
construction given in  \cite{PE3} leads also to strongly chaotic
spaces. Indeed, if $\{ X_{t}: t \in T \}$, $X$ and
$a_{1},a_{2},\ldots$ are as in Proposition 3.3,  both sets $
\bigcup _{j=1}^{\infty}\{a_{2j}\}$ and $ \bigcup
_{j=1}^{\infty}\{a_{2j-1}\}$ are dense in $X$, and
 $s_{1},s_{2}, \ldots $ is a fixed
sequence of different elements of $T$, then for every $t \in T
\setminus \{s_{1},s_{2}, \ldots \}$, the space
$Z_{t}=L(X,K_{i},a_{i})$, where $K_{2j}=X_{s_{j}}$ and
$K_{2j-1}=X_{t}$, is strongly chaotic  and  if $t \neq t'$ then no
open subset of  $Z_{t}$ embeds in $Z_{t'}$.

(B)  Some new examples of strongly chaotic continua were
constructed recently in  \cite{JK}.

\medskip

{\bf  Lemma 3.5.} {\it If $E$ and $A$ are two  continua such that
$E$ is  strongly chaotic and no open nonempty subset of $E$ embeds
in $A$, then every pseudosuspension $M(E,A,a)$ is strongly
chaotic. }
\medskip

{\bf Proof.} Let $p: M(E,A,a)\to E$ be an atomic mapping described
in Remark 2.6. Suppose that there exist  two disjoint subsets $U$
and $V$ of $M(E,A,a)$, with $U$ being nonempty and open, and a
homeomorphism $f$  of $U$ onto $V$. Since $p^{-1}(a)$ is boundary
and closed in $M(E,A,a)$, we can assume that $U$ is a subset of $
p^{-1}(E \setminus \{a\})$. Since no open subset of $E$ embeds in
$A$, then $V$ is not contained in $p^{-1}(a)$.
 Thus $V' = V
\setminus p^{-1}(a)\neq \emptyset $,  hence $f^{-1}(V')=U'$ is a
nonempty open subset of $p^{-1}(E \setminus \{a\})$ homeomorphic
to $V'$ - a contradiction, since  $p^{-1}(E \setminus \{a\})$ is
homeomorphic to  $ E \setminus \{a\}$ and $E$ is strongly chaotic.

\section{Classes of HI strongly chaotic continua without common models}

\medskip

In this section we will prove the results stated in the abstract.

\medskip

{\bf  Theorem 4.1.} {\it Let $\{W_{K}: K \in {\mathcal S}\}$
 be a family from Definition 2.1 and let $n$ be a natural number.  Then there exist a
countable set ${\mathcal S}_{0} \subset {\mathcal S}$ and a family
$\{Y(W_{K}): K \in {\mathcal S} \setminus {\mathcal S}_{0}\}$
consisting of HI strongly chaotic $n$-dimensional Cantor manifolds
such that, for every $K \in {\mathcal S} \setminus {\mathcal
S}_{0} $,  $Y(W_{K})$ admits a mapping onto $W_{K}$. Moreover, if
$n \geq 2$, then every subcontinuum of $Y(W_{K})$ of dimension
greater than 1 admits a   mapping onto $W_{K}$. }

\medskip

{\bf Proof.} Let $E$ be a 1-dimensional HI Cook continuum
satisfying (6) and let $\{ X_{i}: i \in N \}$ be a sequence of
disjoint nontrivial subcontinua of $E$. Then $X_{i}$ does not
embed in $X_{j}$ for any $i \neq j $. For $K \in {\mathcal S}$ let
$X(W_{K})$ be an HI $n$-dimensional Cantor manifold constructed in
Lemma 2.2 and let $a_{1},a_{2},\ldots$ be a dense sequence of
points in $X(W_{K})$ belonging to different composants of
$X(W_{K})$. By (4), there exists a countable set ${\mathcal S}_{0}
\subset {\mathcal S}$ such that for any $K \in {\mathcal S}
\setminus {\mathcal S}_{0}$ and $i \in N$, every mapping from
$X_{i}$ into $X(W_{K})$ is constant (in particular, $X_{i}$ does
not embed into $X(W_{K})$). Fix $K \in {\mathcal S} \setminus
{\mathcal S}_{0}$. Then by Proposition 3.3 and Theorem 2.5 (iv),
(vi),
 the space $Y(W_{K})=L(X(W_{K}),X_{i},a_{i})$ is an HI
strongly chaotic $n$-dimensional Cantor manifold. The space
$Y(W_{K})$ has an atomic mapping $p_{K}$ onto $X(W_{K})$ such that
$p_{K}^{-1}(a_{i})$ is homeomorphic to $X_{i}$ and $p_{K}^{-1}(x)$
is a one-point set for $x \not \in \{a_{i}: i \in N \}$. Since
$X(W_{K})$ maps onto $W_{K}$ then $Y(W_{K})$ can also be mapped
onto $W_{K}$. Now, let $n \geq 2$, and let $Z \subset Y(W_{K})$ be
a subcontinuum of $Y(W_{K})$ of dimension greater than 1. Since
every $p_{K}^{-1}(a_{i})$ is 1-dimensional, then $p_{K}(Z)$ is a
nontrivial subcontinuum of $X(W_{K})$, hence it can be mapped onto
$W_{K}$.

\medskip

{\bf  Theorem 4.2.} {\it Let $\{W_{K}: K \in {\mathcal S}\}$
 be a family from Definition 2.1.  Then there exist a
countable set ${\mathcal S}_{0} \subset {\mathcal S}$ and a family
$\{Y(W_{K}): K \in {\mathcal S} \setminus {\mathcal S}_{0}\}$
consisting of HI strongly chaotic hereditarily strongly
infinite-dimensional Cantor manifolds such that, for every $K \in
{\mathcal S} \setminus {\mathcal S}_{0} $,  $Y(W_{K})$ admits a
  mapping onto $W_{K}$.}

  \medskip

{\bf Proof.}  Let   $\{ X_{i}: i \in N \}$ be a sequence of HI
hereditarily strongly infinite-dimensional Cantor manifolds such
that $X_{i}$ does not embed in $X_{j}$ for any $i \neq j $ (see
Theorem 3.1). For $K \in {\mathcal S}$ let $X(W_{K})$ be an HI
hereditarily strongly infinite-dimensional Cantor manifold
constructed in Lemma 2.3 and let $a_{1},a_{2},\ldots$ be a dense
sequence of points in $X(W_{K})$ belonging to different composants
of $X(W_{K})$.  By (4), there exists a countable set ${\mathcal
S}_{0} \subset {\mathcal S}$ such that  $X_{i}$ cannot be mapped
onto $W_{K}$ for every $K \in {\mathcal S} \setminus {\mathcal
S}_{0}$ and every $i = 1,2,\ldots$. Fix $K \in {\mathcal S}
\setminus {\mathcal S}_{0}$. Since every nontrivial subcontinuum
of $X(W_{K})$ maps onto $W_{K}$, then no $X_{i}$ embeds into
$X(W_{K})$.  Thus, by Proposition 3.3,
$Y(W_{K})=L(X(W_{K}),X_{i},a_{i})$ is strongly chaotic. By Theorem
2.5 (iv), (vii) and (viii), the space $Y(W_{K})$ is an HI
hereditarily strongly infinite-dimensional Cantor manifold. Since
$Y(W_{K})$ has a mapping onto $X(W_{K})$, $Y(W_{K})$ can also be
mapped onto $W_{K}$.

\medskip

{\bf  Theorem 4.3.} {\it Let $\{W_{K}: K \in {\mathcal S}\}$
 be a family from Definition 2.1. For every infinite ordinal number $\alpha
< \omega _{1} $ and for eve\-ry spiral $W_{K}$, where $K \in
{\mathcal S}$, there exists a strongly chaotic HI continuum
$Y(W_{K})$ with ${\rm ind}Y(W_{K}) = \alpha $ (respectively, ${\rm
Ind}Y(W_{K})=\alpha$), which can be mapped onto $W_{K}$. }

\medskip

{\bf Proof.} For a countable infinite ordinal $\alpha$, let $\{
X_{i}: i \in N \}$ be a sequence of continua with small
(respectively, large) transfinite dimension equal to $\alpha$ such
that $X_{i}$ does not embed in $X_{j}$ for any $i \neq j $. Fix $K
\in {\mathcal S}$. Then by Lemma 2.2 there exists a 1-dimensional
HI continuum $X(W_{K})$  which admits a   mapping onto $W_{K}$.
Let $a_{1},a_{2},\ldots$ be a dense sequence of points in
$X(W_{K})$ belonging to different composants of $X(W_{K})$. For
every $i \in N$ the continuum $X_{i}$  is infinite-dimensional, so
it does not embed in $X(W_{K})$. Thus, by Proposition 3.3, the
space $Y(W_{K})=L(X(W_{K}),X_{i},a_{i})$ is  strongly chaotic.
 By Theorem
2.1 (iv) and (v), the space $Y(W_{K})$ is an HI continuum with the
small (respectively, large) transfinite dimension equal to
$\alpha$. Since $Y(W_{K})$ has an atomic mapping onto $X(W_{K})$,
then $Y(W_{K})$ can also be mapped onto $W_{K}$.

\medskip

 From Theorems 4.1 - 4.3 we obtain immediately the following
 corollary.

 \medskip

{\bf  Corollary 4.4. } {\it Let ${\mathcal K}$ be one of the following  classes of spaces:
\begin{lister}
\item[\em 1) ] the class of strongly chaotic $n$-dimensional HI Cantor ma\-ni\-folds, for any natural number $n$,
\item[\em 2) ] the class of strongly chaotic HI hereditarily strongly infinite-dimensional Cantor ma\-ni\-folds,
\item[\em 3) ] the class of strongly chaotic countable-dimensional HI continua of a given transfinite dimension
(small or large) $\alpha < \omega _{1}$.
\end{lister}
Then ${\mathcal K}$ does not have a common model.}

\medskip

{\bf  Remark 4.5.} Observe that, by (1),  for every countable collection $\{X_{i}: i \in N\}$ of continua there exists a spiral
$W_{K}$, where $ K \in {\mathcal S} \setminus {\mathcal S}_{0}$ for any given countable ${\mathcal S}_{0}$, such that no member of
this collection can be mapped onto $W_{K}$. From this fact and Theorems 4.1 - 4.3,\linebreak applying the transfinite induction one
obtains that
\medskip

(7) for every class ${\mathcal K}$ of continua described in
Corollary 4.4, there exists a transfinite sequence $\{Z_{\alpha}:
\alpha < \omega _{1} \}$ of continua from ${\mathcal K}$ such
that, for every $\alpha < \beta < \omega _{1}$, there is no
  mapping from $Z_{\alpha }$ onto $Z_{\beta}$.

\section{Compactifications of the ray as images of HI strongly chaotic continua}

Below we will strengthen some results of the preceding section by
using
   a method suggested by the referee
of the paper \cite{PE1}. This method is based on an idea of
J.T.Rogers, Jr. \cite{Ro}, who applied a theorem of Mazurkiewicz
from \cite{Ma} to show that for every continuum $X$ of dimension
$\geq 2$ and every Waraszkiewicz spiral $W$ there exists a
continuum $X_{W}\subset X$  which admits a mapping $f_{W}$ onto
$W$. The idea of the referee (see \cite{PE1}, Remark 5.2) is that
one can  construct a pseudosuspension $M(E,X_{W},a)$ of $X_{W}$
over any HI continuum $E$ at some point $a$ in such a way that
$f_{W}$ extends to $\tilde{f}_{W}$ from $M(E,X_{W},a)$ onto $W$.
This method, described in detail in the following Lemmas 5.1 and
5.2, allows in particular to obtain HI
 hereditarily strongly
infinite-dimensional Cantor manifolds (or HI $n$-dimensional
Cantor manifolds) which can be mapped onto any given Waraszkiewicz
spiral, or, more generally, onto any given  compactification of
the ray. Recall that, as proved in \cite{Be2} and \cite{AB}, every
continuum is a remainder in some metric compactification of the
ray (cf. the proof of Lemma 5.1 below).

\medskip

{\bf  Lemma 5.1.} {\it  Let $W = L \cup S$ be a compactification
of the ray $L$ with the remainder $S$. Let $X$ and $A$  be two
continua, $a \in X$ and $f:A \to W$ be a   mapping of $A$ onto
$W$.  Suppose that there exists a sequence $A_{1}\subset A_{2}
\subset \ldots$ of subcontinua of $A$ contained in $f^{-1}(L)$
such that the union $\bigcup _{i=1}^{\infty}A_{i}$ is dense in
$A$. Then there exists a pseudosuspension $M(X,A,a)$ which admits
a   mapping $\tilde{f}: M(X,A,a) \to W$ onto $W$.}

\medskip

{\bf Proof.} We adapt the proof of Theorem from the paper
\cite{AB} by J.M.Aarts and P.van Emde Boas.

Suppose that $L$ is an  image of the half-line under a
homeomorphism $h: [0,+\infty ) \to L$, and fix a metric $\sigma $
in $W$. For $z',z'' \in L$, let $L(z',z'')$ denote the arc in $L$
with the end points $z'$ and $z''$.

 One can assume that $A$ is a subset of the
Hilbert cube $I^{\infty}$ with a metric $d$. Choose a countable
dense subset $\{a_{1},a_{2},\ldots \}$ of $\bigcup
_{i=1}^{\infty}A_{i}$. For every $i$, there exists $j(i) $ such
that $a_{i}, a_{i+1} \in A_{j(i)}$.

 There exists $m \in N$ such that $f(A_{j(i)})$ is contained in  $h([0,m])
$. Since $h([0,m])$ is an arc in $L$,  there is an $\epsilon _{i}
> 0$ such that if for $z', z'' \in h([0,m])$ we have $\sigma
(z',z'')< \epsilon _{i}$ then the diameter of the arc $L(z',z'')$
in $h([0,m]) $ is less than $\frac{1}{i}$. Thus, from the uniform
continuity of $f \mid A_{j(i)}$, there is a $\delta _{i} <
\frac{1}{i}$ such that if $d(a',a'')< \delta _{i}$ for $a',a'' \in
A_{j(i)}$ then $\sigma (f(a'),f(a'')) < \epsilon _{i}$ and so the
diameter of $L(f(a'),f(a''))$ is less than $\frac{1}{i}$.

Since $A_{j(i)}$ is a continuum, there is a finite $\delta
_{i}$-chain $a_{i}= a_{i}^{1}, a_{i}^{2}, \ldots , a_{i}^{n_{i}}
=a_{i+1}$ of points of $A_{j(i)}$ from $a_{i}$ to $a_{i+1}$, i.e.,
$d(a_{i}^{j},a_{i}^{j+1})< \delta _{i} < \frac{1}{i}$. By the
choice of $\delta _{i}$,
 the diameter of the arc
$L(f(a_{i}^{j}),f(a_{i}^{j+1}))$ is less than $\frac{1}{i}$ for
every $j = 1,2,\ldots , n_{i}-1$.

 By arranging  these
$\delta _{i}$-chains into a sequence  $a_{1}^{1}, \ldots ,
a_{1}^{n_{1}},a_{2}^{1}, \ldots ,a_{2}^{n_{2}}, a_{3}^{1}, \ldots
$ we  obtain a countable dense subset $B=\{b_{1},b_{2},\ldots\}$
of the set  $\bigcup _{i=1}^{\infty}A_{i}$ with $\lim _{k \to
\infty }d(b_{k},b_{k+1})=0$ such that the diameters of the arcs
$L(f(b_{k}),f(b_{k+1}))$ tend to 0 when $k \to \infty$. Now, let
us choose a sequence of points $\{x_{1},x_{2},\ldots \}$ of $X
\setminus \{a\} $  such that $\rho (a,x_{k+1})< \rho (a,x_{k})$
and $\rho (a,x_{k})< \frac{1}{k}$ for $k = 1,2,\ldots$, where
$\rho$ is a metric in $X$.

As in \cite{AB}, define $g : X \setminus \{a\} \to I^{\infty}$ in
the following way. Put $g(x)=b_{1}$, if $\rho (a,x)\geq \rho
(a,x_{1})$, and $g(x)=(1-t) \cdot b_{k}+tb_{k+1} $, if $\rho
(a,x)=(1-t) \cdot \rho (a,x_{k}) + t\cdot \rho (a,x_{k+1})$, where
$0\leq t \leq 1 $ (note that $g(x_{k})=b_{k}$). Let $G = \{
(x,g(x)): x \in X\setminus \{a\} \}$ be the graph of $g $.

Then the closure $\overline{G}$ of $G$ in $X \times I^{\infty }$
is a compactification of $G$ with the remainder $\{a\}\times A$.
Moreover, the  mapping $p : \overline{G} \to X$, being the
restriction of the projection of $X \times I^{\infty}$ onto $X$,
is an atomic mapping such that $p^{-1}(a)$ is  homeomorphic to $A$
and
   $p \mid p^{-1}(X \setminus \{a\}) : p^{-1}(X \setminus \{a\}) \to X \setminus \{a\} $ is
a homeomorphism.  Thus $M(X,A,a)=\overline{G}$ is a
pseudosuspension of $A$ over $X$ at $a$.

 Let $h_{k}: [0,1]\to
L(f(b_{k}),f(b_{k+1}))$
 be a homeomorphism onto the arc $L(f(b_{k}),f(b_{k+1}))$ such that $h_{k}(0)=f(b_{k})$ and $h_{k}(1)=f(b_{k+1})$.
Let us define  $\tilde{f}:\overline{G}\to W$ in the following way.
Put  $\tilde{f}((a,y))=f(y)$ for $y \in A$. If $x \in X \setminus
\{a\} $, then let
$\tilde{f}((x,g(x)))=\tilde{f}((x,b_{1}))=f(b_{1})$ if $\rho
(a,x)\geq \rho (a,x_{1})$, and $\tilde{f}((x,g(x)))=h_{k}(t)$ if
$\rho (a,x)=(1-t) \cdot \rho (a,x_{k}) + t \cdot \rho (a,x_{k+1})$
(note that $\tilde{f}((x_{k},b_{k}))=f(b_{k})$) . Then $\tilde{f}$
is continuous and onto $W$. To see that $\tilde{f}$ is continuous,
let us check only that if $(y_{n},t_{n})\in G$, where $y_{n} \in X
\setminus \{a\}$ and $t_{n} \in I^{\infty}$, and $(y_{n},t_{n})\to
(a,t)$, where $t \in A$, then $\tilde{f} (y_{n},t_{n}) \to
\tilde{f}(a,t) = f(t)$.

For every $n \in N$  let $k_{n}$ be such that $\rho
(x_{k_{n}+1},a) < \rho (y_{n},a) \leq \rho (x_{k_{n}},a) $. Then
$\lim (x_{k_{n}},b_{k_{n}})= \lim (y_{n},t_{n})= (a,t)$, since
$t_{n}$ belongs to the interval with the end points $b_{k_{n}}$
and $b_{k_{n}+1}$ and $d(b_{k},b_{k+1})\to 0$, so $\lim
d(t_{n},b_{k_{n}})=0$. From the definition of $\tilde{f}$,
$\tilde{f} (y_{n},t_{n})\in L(f(b_{k_{n}}),f(b_{k_{n}+1}))$. Since
the diameters of $L(f(b_{k_{n}}),f(b_{k_{n}+1}))$ tend to 0, then
$\lim \tilde{f} (y_{n},t_{n})=\lim f(b_{k_{n}})$. Since $\lim
b_{k_{n}}=t$, then by the continuity of $f$, $\lim
f(b_{k_{n}})=f(t)=\tilde{f}(a,t)$, which ends the proof.

\medskip

 Note that the mapping
$\tilde{f}:M(X,A,a)\to W$ constructed in Lemma 5.1 is in fact an
extension of $f : A \to W$.

\medskip

{\bf  Lemma 5.2.} {\it Let $W = L \cup S$ be a compactification of
the ray $L$ with a remainder $S$. Let $f:Y \to W $ be a mapping
from a conti\-nuum $Y$ onto   $W $. Then there exists a sequence
$A_{1}\subset A_{2} \subset \ldots$ of subcontinua of $Y$
contained in $f^{-1}(L)$ such that if $A$ is the closure of the
union $\bigcup _{i=1}^{\infty}A_{i}$, then $f \mid A$ is onto
$W_{K}$. }

\medskip

{\bf Proof.} Suppose that $L$ is the image of the half-line under
a homeomorphism $h: [0,+\infty ) \to L$. Let $I_{n}= h([0,n))$ for
$n \in N$.
 Inductively
we will construct continua $A_{0} \subset A_{1}\subset A_{2} \subset \ldots$ and points $a_{i} \in A_{i}$ such that $f(a_{i})=h(i)$
for every $i = 0,1,2,\ldots$ Let $a_{0}$ be a point such that $f(a_{0}) = h(0)$ and $A_{0}=\{a_{0}\}$. Suppose that such continua
$A_{i}$ and points $a_{i}$ are already constructed for $i \leq n$. Consider the open  set $U_{n+1}= f^{-1}(I_{n+1}) $ and let
$A_{n+1}$ be a component of the set $\overline{U_{n+1}}$ which contains $a_{n}$. Then $A_{n+1} $ intersects the boundary of $U_{n+1}$,
by a theorem of Janiszewski (see \cite{Ku}, \S 47, III, Theorem 1). Let $a_{n+1}\in A_{n+1} \cap Fr(U_{n+1})$. We have then
$f(a_{n+1})=h(n+1)$. Also, $f(A_{i})$ is a continuum in $W$ containing $h(0)$ and $h(i)$, hence $f(A_{i})\supset I_{i}$. Now, if $A =
\overline{\bigcup _{i=1}^{\infty}A_{i}}$, then the map $f \mid A : A \to W$ is onto.

\medskip

Recall that a mapping is {\it monotone} if each of its fibers is a
continuum. A mapping $f : X \to Y$ is {\it weakly confluent} if
for each continuum $B\subset Y$ there exists a continuum $A\subset
X$ such that $f(A)=B$.

\medskip

{\bf Proposition 5.3.} {\it For every continuum $W$ and every HI continuum $Z$ of dimension~$\geq\hspace{-.1ex} 2$, there
exists a subcontinuum $Y\subset Z$ which can be mapped onto $W$.}

\medskip

{\bf Proof.} By a theorem of Hurewicz, there exists a monotone
mapping $h$ from some metric curve $M$ onto the Hilbert cube $Q$
(for the proof see \cite{MT}, (19.1)). We can assume that $M
\subset I^{3}$ and $W \subset Q$; then the preimage $h^{-1}(W)$ is
a continuum in $I^{3}$.

By a theorem of Kelley \cite{Ke}, there is a monotone open mapping
$f$ of $Z$ onto an infinite-dimensional HI continuum $K$. Since
dim$K\geq 3$, by a theorem of Mazurkiewicz \cite{Ma} there is a
weakly confluent mapping $g$ of $K$ onto $I^{3}$. Since the
composition $gf$ is weakly confluent, there exists a continuum
$Y\subset Z$ such that $gf(Y)=h^{-1}(W)$. The composition $hgf$
maps $Y$ onto $W$.



\medskip

{\bf Remark 5.4.} (A) Proposition 5.3 can be strengthened to the
following effect: {\it for every continuum $W$ and every HI
continuum $Z$ of dimension $\geq 2$, there exists a subcontinuum
$Y\subset Z$ which admits a weakly confluent mapping  onto $W$}.
Indeed, by a theorem of Ma\'{c}\-ko\-wiak and Tymchatyn \cite{MT},
(19.3), there is an HI curve $L$ which admits a weakly confluent
map $g$ onto $W$. Then, by Proposition 5.3, there are a
subcontinuum $Y\subset Z$ and a map $f$ from $Y$ onto $L$. Since
$L$ is HI, a theorem of Cook \cite{Co} implies that $f$ is weakly
confluent. Thus, the composition $gf$ is also weakly confluent.

One can give also  a more direct proof of the above fact, similar
to the proof  of (19.3) in \cite{MT}.

(B) Let us show how Proposition 5.3 can be used in the proof of
the following theorem of  J.Kra\-sin\-kie\-wicz \cite{Kra}:  {\it
the hyperspace $C(X)$ of subcontinua of a nondegenerate HI plane
continuum $X$ is two-dimensional} (cf. also \cite{Kra2}, p.42,
where a similar argument is used).

First we prove that  {\it no HI continuum of dimension $\geq 2$ is an image of a plane continuum.} Suppose $f$ is a map from a plane
continuum $X$ onto an HI continuum $Z$ of dimension~$\geq 2$. Then, by Proposition 5.3, $Z$ contains a continuum $Y$ which admits a
map  onto a diadic solenoid. A theorem of Cook \cite{Co} implies that $f$ is weakly confluent, and hence, $X$~contains a continuum $A$
such that $f(A)=Y$. A contradiction, since a diadic solenoid is not an image of a plane continuum (M.K.Fort, Jr.\ \cite{Fo}, cf.\ also
W.T.Ingram~\cite{In2}).

Now, it suffices to use a theorem by C.Eberhart and S.B.Nadler, Jr.\ \cite{EN} who proved that, if $X$ is a nondegenerate HI
continuum, then (i) \mbox{$\dim C(X)=2$} or $\dim C(X)=\infty$, and (ii) $\dim C(X)=2$ if and only if each monotone open
image of $X$ is one-dimensional (note that a monotone open image of an HI continuum is HI).

\medskip

{\bf  Theorem 5.5.} {\it For every compactification $W$ of the ray
(hence, for every Waraszkiewicz spiral), there exists a strongly
chaotic  HI hereditarily strongly infinite-dimensional Cantor
 manifold $Y(W)$ which can be mapped onto $W$.}

 \medskip

 {\bf Proof.} By Theorem 3.1 there exist two strongly chaotic HI
  hereditarily strongly infinite-dimensional Cantor
 manifolds $E$ and $Z$ such that no open
 nonempty subset of $E$ embeds in $Z$. By Proposition 5.3 there is a continuum $Y \subset
 Z$ and a mapping $f$ of $Y$ onto $W$. By Lemma 5.2 there exists a continuum $A \subset
 Y$ sa\-ti\-sfying the assumptions of Lemma 5.1, hence there exists a
 pseudosuspension $M(E,A,a)$ and a   mapping $\tilde{f}$
 of $M(E,A,a)$ onto $W$. Since $A \subset Z$, no open nonempty subset
 of $E $ embeds in $A$, so by Lemma 3.5 the space $M(E,A,a)$ is
 strongly chaotic. By Remark 2.6 (see (iv)', (vii)' and (viii)'), $M(E,A,a)$
 is an HI hereditarily strongly infinite-dimensional Cantor manifold, so we can put $Y(W)=M(E,A,a)$.

\medskip

Similarly, one  obtains the following  theorem.

\medskip

{\bf Theorem 5.6.} {\it Let $n\geq 2$ be a fixed natural number. Then for every compactification $W$ of the ray (hence, for every
Waraszkiewicz spiral), there exists a strongly chaotic HI $n$-dimensional Cantor manifold $Y(W)$ which can be mapped onto $W$.}

\medskip

 {\bf Proof.} The proof is similar to the proof of Theorem 5.5. By Theorem 3.1 there exist two strongly chaotic HI
 $n$-dimensional Cantor manifolds $E$ and $Z$ such that no open
 nonempty
 subset of $E$ embeds in $Z$. By Proposition 5.3 there is a continuum $Y \subset
 Z$ and a mapping $f$ of $Y$ onto $W$. By Lemma 5.2 there exists a continuum $A \subset
 Y$ sa\-ti\-sfying the assumptions of Lemma 5.1, hence there exists a
 pseudosuspension $M(E,A,a)$, where $a$ is an arbitrary point of $E$, and a   mapping $\tilde{f}$
 of $M(E,A,a)$ onto $W$. Since $A \subset Z$, no open nonempty subset
 of $E $ embeds in $A$, so by Lemma 3.5 the space $M(E,A,a)$ is
 strongly chaotic. By Remark 2.6 (see (iv)' and (vi)'), $M(E,A,a)$ is an HI $n$-dimensional Cantor
 manifold, so we can put $Y(W)=M(E,A,a)$.

\medskip

{\bf  Theorem 5.7.} {\it Let $\alpha < \omega _{1} $ be an
infinite ordinal number. Then for every compactification $W$ of
the ray (hence, for every Waraszkiewicz spiral) there exists a
strongly chaotic HI continuum $Y(W)$ with ${\rm ind}Y(W) = \alpha
$ (respectively, ${\rm Ind}Y(W)=\alpha$), which can be mapped onto
$W$. }

\medskip

{\bf Proof.} From Theorem 3.1 it follows that there exists an HI
strongly chaotic continuum $E$ with ind$E=\alpha$ (Ind$E =
\alpha$) every nonempty open subset of which is
infinite-dimensional. From the proof of this theorem (cf.
Proposition 3.3) it follows that $E$ is of the form $X_{\sigma}
=L(X,X_{s_{i}},a_{i})$, where ind$X=1$ and ind$X_{s_{i}}=\alpha$
(Ind$X_{s_{i}}=\alpha$) for every $i$. Let $p:L(X,X_{s_{i}},a_{i})
\to X$ be an  atomic mapping described in Theorem 2.5. Let $Z$ be
any three-dimensional HI continuum. As in the proof of Theorem
5.5, applying Proposition 5.3 and Lemma 5.2, one obtains a
continuum $A \subset Z$ that satisfies the assumptions of Lemma
5.1, where we take $X = E = X_{\sigma}$. Choose $a \in X \setminus
\bigcup _{i=1}^{\infty}\{a_{i}\}$ and observe that $p^{-1}(a)$ is
a singleton. Then by Lemma 5.1 there exists a pseudosuspension
$M(E,A,p^{-1}(a))$ and a   mapping $\tilde{f}$
 of $M(E,A,p^{-1}(a))$ onto $W$. The continuum $M(E,A,p^{-1}(a))$ is HI, since $E$ and
 $A $ are HI (see Remark 2.6 (iv)'). Note that $M(E,A,p^{-1}(a))$ admits an atomic mapping
 $g $ onto the 1-dimensional continuum $X$ such that $g^{-1}(a)$ is a
 copy of $A$ and $p^{-1}(a_{i})$ is a copy of $X_{s_{i}}$ with
 ind$X_{s_{i}}=\alpha$ (Ind$X_{s_{i}}=\alpha$). It follows that
 ind$M(E,A,p^{-1}(a))=\alpha$ (Ind$M(E,A,p^{-1}(a))=\alpha$) (see the proof of Theorem 3.2 (v) in
 \cite{PR}). Recall that  $E$ is strongly chaotic and no open nonempty subset $U$ of E
embeds into the continuum $A$, since dim$U = \infty$ and dim$A\leq
3$. Thus by Lemma 3.5 the space $Y(W_{K})=M(E,A,p^{-1}(a))$ is
strongly chaotic.

\medskip

{\bf Acknowledgements.} We would like to thank the anonymous
referee of the paper \cite{PE1} whose idea has made an important
contribution to this paper.

\end{document}